\documentclass{amsart}
\usepackage{tabularx}
\usepackage{enumerate}
\usepackage{enumitem,amsmath,amsfonts,amssymb,xcolor}
\usepackage{latexsym,amsfonts,euscript,amssymb,gensymb}
\usepackage[numbers,sort&compress]{natbib}
\usepackage[flushmargin]{footmisc}
\newtheorem{lemma}{\bf Lemma}[section]

\newtheorem{defi}[lemma]{\bf Definition}
\newtheorem{prop}[lemma]{\bf Proposition}
\newtheorem{thm}[lemma]{\bf Theorem}

\newtheorem{exa}[lemma]{\bf Example}
\newtheorem{cor}[lemma]{\bf Corollary}

\newtheorem{con}[lemma]{\bf Conjecture}
\newtheorem{pro}[lemma]{\bf Problem}
\newtheorem{property}[lemma]{\bf Property}

\newcommand{\SL}{{\operatorname{SL}}}
\newcommand{\PGL}{{\operatorname{PGL}}}

\newcommand{\PSL}{{\operatorname{PSL}}}

\newcommand{\PSU}{{\operatorname{PSU}}}

\DeclareMathOperator{\Aut}{Aut}

\DeclareMathOperator{\EPPO}{EPPO}
\DeclareMathOperator{\EPO}{EPO}
\DeclareMathOperator{\diam}{diam}
\DeclareMathOperator{\COE}{COE}
\DeclareMathOperator{\CLT}{CLT}
\DeclareMathOperator{\OC}{OC}

\input xy
\xyoption{all}

\title[]{Quantitative characterization of finite simple groups: a complement}

\author{Wujie Shi}
\address{\rm{Wujie Shi. School of Mathematics and Big Date, Chongqing University of Arts and Sciences, Chongqing, P.R.China,  402160;
    School of Mathematical Sciences, Soochow University, Suzhou, Jiangsu, P.R.China, 215006.}}
\email{wjshi@suda.edu.cn; shiwujie@outlook.com}

\thanks{Project supported by the NSF of China (Grant No.11171364 and 11271301).}

\date{}


\begin{document}

\setlength{\parskip}{2mm}

\maketitle

\begin{abstract}
In this paper, we summarize the work on the characterization of finite simple groups and the study on finite groups with ``the set of element orders" and ``two orders" (the order of group and the set of element orders). Some related topics, and the applications together with their generalizations are also discussed. The original version of this article was published in Chinese in the journal Scientia Sinica Mathematica, no.53(2023), pp.931-952. This revised and expanded version has corrected several errors and added quite a few contents. Especially, it is pointed that this work has applications in mathematics and computational complexity theory.
\end{abstract}

\section{Introduction}
Group theory is an important branch of mathematics, and it has important applications in mathematics, physics, chemistry and other fields. The classification theorem for the finite simple groups, which was completed in 2004, is one of the most important mathematical achievements of the 20th century. In terms of the number of participating researchers, the number of published papers, and the total length of the proof (more than 15,000 pages), it is unprecedented in the history of mathematics. The first paper classifying an infinite family of finite simple groups, starting from a hypothesis on the structure of certain proper subgroups, was published by W. Burnside in 1899.

In the history of group theory, Burnside conjectured that every nonabelian simple group has even order in 1911. R. Brauer suggested classifying simple groups by using their centralizers of involutions in 1954. Brauer and K. A. Fowler
proved that if $G$ is a simple group of even order and $\tau$ an involution in $G$, then $|G|\leqslant (|C_G(\tau)|^2)!$. Therefore, given $C_G(\tau)$, there are only finitely many simple groups with centralizers of involutions isomorphic to $C_G(\tau)$, and this is Brauer program. Brauer-Fowler theorem proved that if the centralizer of an involution in simple group is known, then the simple group can be determined (See \cite{22}). It is well-known that a group of odd order does not have element of order $2$. Therefore, in order to use Brauer program, we must first prove that the orders of  nonabelian simple groups are even. This is equivalent to proving that groups with odd orders are solvable, which is precisely what W. Feit and J. G. Thompson had later proved.

In 1957, M. Suzuki started to prove Burnside's conjecture. He studied CA groups, which are groups with the centralizer of every nontrivial element is abelian. He published a groundbreaking paper in which he proved that all CA groups of odd order are solvable. Subsequently, he classified all the simple CA groups, and more generally, classified the simple groups with centralizers of involutions having a normal Sylow $2$-subgroup. In this process, a family of simple groups of Lie type was discovered, known as the Suzuki simple groups. In 1960, Feit,  Thompson, and M. Hall extended Suzuki's work to CN groups, which are groups with the centralizer of every nontrivial element is nilpotent. They proved that every CN group of odd order is solvable.

To prove the theorem that odd order groups are solvable, Feit and Thompson first proved that there does not exist a nonabelian simple group of odd order in which every subgroup is solvable. This implies that the minimal counterexample to the solvability of odd order groups is a simple group in which every subgroup is solvable. Although the general approach of its proof is similar to the CA and CN theorems, the details are much more complex, with their final paper spanning 254 pages.

\begin{thm}[Feit-Thompson Theorem (1963)]
Let $G$ be a finite group. If the order $|G|$ of $G$ is odd, then $G$ is solvable.
\end{thm}

The Odd Order Theorem tells us that the quantities of a group can determine its properties.

Recently, many significant topics in group representation theory have been successfully reduced to the cases of simple groups and quasi-simple groups, such as McKay Conjecture and Alperin Conjecture, and so on (See \cite{146}), we urgently need to understand and know more properties of simple groups. The famous researcher R. Solomon in the field of group theory pointed out that: ``The ongoing project to publish a series of more than $12$ volumes presenting a complete proof of this theorem is expected to be completed by 2025." (See the email sent by Solomon to the author on November 2nd, 2020 or \cite {186}).

The Odd Order Theorem provides us with insight that ``In a group, quantities can determine the structure of the group". But in finite groups, not only the ``order of the group", but also the ``element orders" is a basic quantity. Thus, the author first provided a characterization of $A_5$: a finite group $G$ is isomorphic to $A_5$ if and only if the set of element orders in $G$, denoted by $\pi_e(G)$, is $\{1, 2, 3, 5\}$. To generalize this result, it is necessary to add the condition of the ``group order". In 1987, after characterizing all sporadic simple groups using the ``group order" and the ``set of element orders", the author proposed the following conjecture:

\begin{con}[See \cite{165}]\label{1.1}
  Let $G$ be a finite group and $S$ a finite simple group. Then $G \cong S$
if and only if $(1)$ $\pi_e(G)= \pi_e(S)$; $(2)$ $|G|=|S|$, in other words, every finite simple group can be characterized by using only the order of the group and the orders of its elements \textup{(}briefly, ``two orders"\textup{)}.
\end{con}

In 1987, the author posed the above conjecture to Thompson and received his encouragement and full affirmation. In his response, Thompson remarked: ``Good luck with your conjecture about simple groups. I hope you continue to work on it", ``I like your arguments", ``This would certainly be a nice theorem". At the same time, in two letters from Thompson dated 1987 and 1988 respectively, he posed the following problem and conjecture (see reference \cite [Problem 12.37-12.39]{100}).

\begin{defi}\label{1.2}
    Let $G$ be a finite group, let $d$ be a positive integer and write $G(d)=\{x \in G | x^d = 1 \}$. We say that two groups $G_1$ and $G_2$ are of the same order type if and only if $|G_1(d)|=|G_2(d)|$ where $d=1$, $2, \cdots $.
\end{defi}

\noindent\textbf{Thompson Problem (1987)}
\footnote {A private letter from J. G. Thompson dated April 22, 1987 to the author raised this problem, which is also introduced in E. I. Khukhro, V. D. Mazurov.\textit{ Unsolved problems in group theory.} The Kourovka notebook, No. 20. 2022, 12.37, p.58.}
Suppose $G_1$ and $G_2$ are groups of the same order type. Suppose also that $G_1$ is solvable. Is it true that $G_2$ is also necessarily solvable?

In Thompson's letter he pointed out that:``The problem arose initially in the study of algebraic number fields, and is of considerable interest."

Let $G$ be a finite group, set $N(G)=\{n\in\mathbb{Z}^+ \mid G \text{ has conjugacy class } C \text{ with } |C| = n \}$, that is, it is the set of all lengths of the conjugacy classes of $G$.

\noindent\textbf{Thompson Conjecture (1988)}
\footnote{A private letter from J. G. Thompson dated January 4, 1988 to the author posed this conjecture, which is also introduced in E. I. Khukhro, V. D. Mazurov. \textit{Unsolved problems in group theory.} The Kourovka notebook, No. 20. 2022, 12.38, p.59.}
Let $G$ and $M$ be two finite groups and $N(G)=N(M)$, if $M$ is a non-abelian simple group and the centre of $G$ is trivial, then $G$ and $M$ are isomorphic.

This article provides an overview of the proof for Conjecture \ref{1.1}, Thompson Conjecture, the current status of solutions for Thompson Problem, as well as some problems related to quantities and groups, it involves their applications and generalizations, and also includes some unsolved problems.

In fact, Conjecture \ref{1.1} originated from the author's master's thesis. The thesis studied finite groups whose elements with prime order except the identity element (element prime order groups, $\EPO$-groups) (See \cite{184}), and finite groups whose elements with prime power orders (element prime power order groups, $\EPPO$-groups) (see \cite{185}). Moreover, the following theorem is obtained (See \cite{163}):

\begin{thm}\label{A5}\label{1.3}
    Let $G$ be a finite group. Then $G\cong A_5$ if and only if $\pi_e(G)=\{1, 2, 3, 5\}$.
\end{thm}

After learning about this result, Xuefu Duan wrote in a letter to Zhongmu Chen: ``Very interesting, I will also think about it when I have time." The encouragement from both Duan and Zhongmu Chen prompted the author to
characterize more simple groups by using ``the set of element orders".
However, when dealing with simple groups (such as $A_6$) that could not be characterized  only by the ``set of element orders", the natural thought is to add the ``order of group." Consequently, the conjecture \ref{1.1} was posed.

Article \cite{184} was published in Chinese. A similar article was published five years later in \cite{42}, but the main theorem of this article contains some mistakes. In response, the author, together with coauthors, published article \cite{39} as a corrigendum to article \cite{42}. Initially published in Chinese as well, article \cite{185} gained recent attention, prompting the authors of \cite{185} to translate it into English and submited it on arXiv. \cite{149} generalized \cite{184} to investigate finite groups in which all elements, except those in a normal subgroup, have prime order. The recently published paper \cite{LML} is a generalization of the reference \cite{185}.

\section{The spectrum characterization of finite simple groups}

Before characterizing the alternating group $A_5$ only using the set of element orders, the author had used the number $|\pi(G)|$ of prime factors of $|G|$ and impose some restrictive conditions on ``the set of element orders" to characterize simple group $\PSL_2(7)$ (See \cite{160}) and some other simple groups (See \cite{159,161,162,164}). Moreover, it is not difficult to deduce that the aforementioned constraints hold from the ``set of element orders". For instance, the author \cite{160} proved the following result:

\begin{thm}\label{2.1}
    Let $G$ be a finite group satisfying the following conditions:
    \begin{itemize}
        \item [\rm{(1)}] $|G|$ contains at least three different prime factors, that is, $|\pi(G)|\geqslant 3$;
        \item [\rm{(2)}] the order of every nontrivial element in $G$ is either a power of $2$ or a prime different from $5$.
    \end{itemize}
Then G is isomorphic to $\PSL_2(7)$.
\end{thm}

It is easy to deduce that the above conditions $(1)$ and $(2)$ hold from $\pi_e(G)=\{1, 2, 3, 4, 7\}$ ($=\pi_e(\PSL_2(7))$).

The alternating group $A_5$ was characterized using the ``set of element orders" as a condition (See \cite{163}), and the proof only used elementary group theory concepts.
Subsequently, in reference \cite{153}, is also used an elementary method to prove that a finite group $G$ is isomorphic to $A_5$ if and only if $\pi_e(G)=\{1, p, q, r\}$, where $p$, $q$, and $r$ are distinct prime numbers.

Due to the significance of the concept of ``the set of element orders in a group", it is referred to as the ``spectrum" in some subsequent articles.

V. D. Mazurov (See \cite{129}) point out: ``This result opened a wide way for investigations of recognizability of groups."

Obviously, the set of element orders in a group $G$, denoted by $\pi_e(G)$, is a subset of the set of positive integers $\mathbb{Z}^+$. However, conversely, the following question is quite difficult:

\begin{pro}\label{2.2}
    What kind of sets of numbers can become a set of element orders $\pi_e(G)$ of a finite group $G$?
\end{pro}

Zhigang Wang, A. V. Vasil'ev, M. A. Grechkoseeva, and A. Kh. Zhurtov proved the following theorem concerning the spectrum $\pi_e(G)$, which provides a helpful criterion of nonsolvability of a finite group (See \cite[Theorem 1] {WZH}).

\begin{thm}\label{spectrum}
    Let $G$ be a finite group. Suppose that there is a subset $\sigma(G)$ of $\pi(G)$, which contains at least four elements and satisfies the following:
    \begin{itemize}
        \item [\rm{(1)}] $pq\in \pi_e(G)$ for all distinct $p,q\in \sigma(G)$;
        \item [\rm{(2)}] $pqr \not\in \pi_e(G)$ for all pairwise distinct $p, q, r\in \sigma(G)$.
    \end{itemize}
    Then $G$ is nonsolvable. In particular, there are no groups of odd order satisfying these conditions.
\end{thm}

Let $\lambda(G)$ denote the maximum number of pairwise nonadjacent vertices in prime graph GK(G) (See Definition \ref{2.11}). It is obvious that a direct square $L\times L$ satisfies those two conditions of Theorem \ref{spectrum} if the group $L$ satisfies $\lambda(L)\geqslant 4$. Thus the following result can be obtained by the above Theorem \ref{spectrum}
(See \cite [Corollary 2]{WZH}).

\begin{cor}\label{nonsolvable1}
    Let $L$ be a finite group and $\lambda(L)\geqslant 4$. If $G$ is a finite group such that $\pi_e(G)=\pi_e(L\times L)$, then $G$ is nonsolvable.
\end{cor}

In October 1996, the author gave the following definition in a presentation of ``Group Theory Seminar" held at Princeton University: For a subset $\Gamma$ of the set of positive integers, the function $h$ can be defined as follows, $h(\Gamma)$ is the number of isomorphism classes of groups $G$ such that $\pi_e(G)= \Gamma$.

\begin{defi}[See \cite{175}]\label{2.3}
     For a given group $G$, we have $h(\pi_e(G))\geqslant 1$. A group $G$ is called recognizable if $h(\pi_e(G))=1$. A group $G$ is called almost recognizable \textup{(}or $k$-recognizable\textup{)} if $h(\pi_e(G))=k$ is finite; otherwise $G$ is called unrecognizable.
\end{defi}

 Solvable groups are unrecognizable (see \cite [Theorem 4]{175}). In particular, the following result was obtained (See \cite{67, 125, 136, 170, 177}):

 \begin{thm}\label{2.4}
  If $G$ has a solvable minimal normal subgroup or $G$ is isomorphic to one of the following groups: $A_6$, $A_{10}$, $L_3(3)$, $U_3(3)$, $U_3(5)$, $U_3(7)$, $U_4(2)$, $U_5(2)$, $J_2$, $S_4(q)$ where $q\neq 3^{2m+1}$ and $m>0$, then $G$ is unrecognizable.
 \end{thm}

In April 2010, at the Ischia Group Theory Conference held in Italy, the author presented the following theorem for recognizable groups (See \cite{177}):

\begin{thm}\label{2.5}
 Let $G$ be one of the following simple groups:
 \begin{itemize}
     \item [\rm{(1)}] Alternating group $A_n$ where $n\neq 6$, $10$.
     \item [\rm{(2)}] Sporadic simple group $S$ where $S\neq J_2$.
     \item [\rm{(3)}] Simple group of Lie type:
     \begin{itemize}
         \item $L_2(q)$ with $q\neq 9$; $L_3(2^m)$ for $m\geqslant 1$; $L_3(q)$ where
               $3<q\equiv 2(\bmod~5)$ and $(6,\frac{q-1}{2})=1$;
         \item $U_3(2^m)$ where $m\geqslant 2$; $L_4(2^m)$ where $m\geqslant 1$;
               $U_4(2^m)$ where $m\geqslant 2$;
         \item $Sz(2^{2m+1})$ with $m\geqslant 1$; $R(3^{2m+1})$ with $m\geqslant 1$;
               ${}^2F_4(2^{2m+1})$ where $m\geqslant 1$; $S_4(3^{2m+1})$ where $m\geqslant 0$;
         \item $B_p(3)$ where $p>3$ is odd prime; $C_p(3)$ where $p$ is odd prime; $D_p(5)$ where $p$ is odd prime;
         \item $D_n(q)$, where $q=2$, $3$ or $5$ for some $n$; $G_2(3^m)$;
         \item ${}^2F_4(2)'$, $L_3(7)$, $L_4(3)$; $L_n(2)$ where $n\geqslant 3$;
               $L_5(3)$, $U_3(9)$, $U_3(11)$, $U_4(3)$, $U_6(2)$, $G_2(3)$;
         \item $G_2(4)$, $S_6(3)$, $O^-_8(2)$, $O^-_{10}(2)$, $F_4(2)$,
               ${}^2E_6(2)$.
    \end{itemize}
 \end{itemize}
 Then $G$ is recognizable.
\end{thm}

\begin{cor}\label{Monster}
 Let $G$ be a finite group and $M$ a Monster group, then $G\cong M$ if and only if  $\pi_e(G)=\{1, 2, \cdots, 36, 38, 39, 40, 41, 42, 44, 45, 46, 47, 48, 50, 51, 52, 54, 55, 56,\\ 57, 59, 60, 62, 66, 68, 69, 70, 71, 78, 84, 87, 88, 92, 93, 94, 95, 104, 105, 110, 119\}$.
\end{cor}

The perface of the book \textit{Moonshine, the Monster and Related Topics} \cite{DCY} points out: ``One of the great legacies of the classification of the finite simple groups is the existence of the Monster. It was the study of this group that first suggested that there might be interesting relations between finite groups and certain elliptic modular functions, and it was this possibility-fuelled by the Conway-Norton conjectures-that led to what one might call the first version of ``Moonshine", that is, the study of class functions on groups with values in a ring of modular functions."

Does Corollary \ref{Monster} make sense for the related research on Monster and related topics?

For the case of $h(\pi_e(G))=2$, Buturlakin (See \cite {230} ) proved that if $G$ and $H$ are finite simple groups such that $\pi_e(H)=\pi_e(G)$ and $H \ncong G$, then $\{G, H\} = \{S_6(2), O^+_8(2)\}$ or $\{G, H\} = \{O_7(3), O^+_8(3)\}$. Moreover, the author of \cite {230} proved that there are no the case of $h(\pi_e(G))=3$ for three different finite simple groups $G$.

Now, the research on the recognizable groups mentioned above has already been summarized comprehensively, see \cite[Table 1-9]{67}, recognizable simple groups are those finite simple groups in which the function $h$ equal $1$ in these tables. Moreover, $h(G)$ is known for all finite simple groups $G$ except for some classical groups in dimensions from $5$ to $36$ over fields of odd characteristic. Furthermore, the following theorem can be found in \cite [Table 10]{67}.

\begin{thm}
   Let $L$ be one of the following groups: $A_6$, $A_{10}$, $J_2$, ${}^3D_4(2)$, $L_3(3)$, $U_3(3)$, $U_3(5)$, $U_3(q)$ where $q\geqslant 7$ and $q$ is special Mersenne, $U_5(2)$, $S_4(3)$, $S_4(2^m)$, $S_4(3^m)$ with $m$ even, $S_4(q)$ where $q=p^m$ and $p\geqslant 5$, $S_8(q)$ for $q=p^m$ and $p\neq 2,7$, $O_9(q)$ where $q=p^m$ Then $L$ is an unrecognizable simple group.
\end{thm}

In recent years, some researchers have characterized the direct square of finite simple groups by spectrum and obtained some interesting results. We say that groups $G$ and $H$ are isospectral if $\pi_e(G)=\pi_e(H)$. Wang, Vasil'ev, Grechkoseeva, and Zhurtov gave the recognizability of the direct square of the Suzuki groups by Corollary \ref{nonsolvable1} and proved the following two theorems (See \cite [Theorems 3, 4]{WZH}):

\begin{thm}
If $q\geqslant 8$ and $q\neq 32$, then the group $Sz(q)\times Sz(q)$ is recognizable by spectrum.
\end{thm}

\begin{thm}
    Let $L=Sz(32)$. Then finite groups isospectral to $L\times L$ are exactly groups of the form $(L\times L)\rtimes \langle \varphi\rangle$, where either $\varphi=1$, or $\varphi$ normalizers each direct factor and induces a field automorphism of order $5$ on it. In particular, up to isomorphism, there are four finite groups isospectral to $L\times L$.
\end{thm}

We know that a nonsolvable simple group with abelian Sylow $2$-subgroups is either a linear group $\PSL_2(q)$, a small Ree group $R(q)={}^2G_2(q)$ where $q=3^{a}$ with $a\geqslant 3$ odd, or the sporadic Janko group $J_1$. R. Brandl and Wujie Shi \cite{21} provided that every finite nonabelian simple group with abelian Sylow $2$-subgroups is recognizable.
The authors \cite{LIT} considered the recognizability for the direct square of the finite nonsolvable simple groups with abelian Sylow $2$-subgroups and established the following result (See \cite [Corollary 1.3]{LIT}):

\begin{thm}
  Suppose that $L$ is a simple group with abelian Sylow $2$-subgroups. Then the direct square $L\times L$ of $L$ is recognizable by spectrum if $L$ is either a small Ree group $R(q)$ or the sporadic Janko group $J_1$, and $L\times L$ is unrecognizable otherwise.
\end{thm}

Mazurov et al. relaxed the condition mentioned above by removing the restriction that the group is ``finite", proving a series of results (See \cite{70,85,114, 115,116,117,118,119,126,127,128,130,131,132,133,135,208,224,225,226,227}). A typical result is as follows (See \cite{113}):

\begin{thm}\label{2.6}
     If the spectrum of a group $G$ is equal to $\{1,2,3,4,7\}$ then $G\cong \PSL_2(7)$.
\end{thm}

Note that if the group $G$ is assumed to be finite, then the conclusion of the above theorem appeared in the literature \cite {160} in 1984. Without the ``finite" condition, a paper with the same conclusion was published in 2007, after a span of 23 years.

For the function $h$, it was conjectured in \cite {148} and \cite[Problem 12.84]{100} that $h(\Gamma) \in \{0,1,\infty\}$.
The authors \cite{40,122,123,138,173,183,192} subsequently provided counterexamples to this conjecture and gave some examples with function $h$ equal to 2. Furthermore, the following question was put forward in \cite{173}: Let $\Gamma$ be an arbitrary subset of the natural numbers, does there exist a positive integer $k$ such that $h(\Gamma) \in \{0, 1, 2,\cdots, k, \infty\}$? If the answer is positive, what is the value of such $k$?
A. V. Zavarnitsine \cite{214} gave an example as follows: For any $r>0$, $h(\pi_e(L_3(7^n)))=r+1$ where $n=3^r$, and these $r+1$ groups are all of the form $L_3(7^n)\langle \rho\rangle$, $\rho$ is a field automorphism of $L_3(7^n)$ ($n= 3^r$), $k= 0, 1, 2,\cdots r$.
It is an example of function $h$ such that for any positive integer $k$, there exist $k$-recognizable groups.
Observing example of group for which the function $h$ is equal to $2$ (See \cite[Theorem 5.4]{170}), the author has raised the question of whether, for any $k$, there exist two finite $k$-recognizable isopectral groups that are section-free, in other words, they are not simple sections of each other (see \cite[Problem 5.2]{177}). Grechkoseeva \cite{66}  pointed out the existence of such groups, which are $G_1 = L_{15}(2^{60}).3$ and $G_2 = L_{15}(2^{60}).5$ (See \cite [Problem A:16.106]{100}). Note that although $G_1$ and $G_2$ mentioned above are section-free, both of them have same normal subgroup $L_{15}(2^{60})$.

The spectrum of a finite group is a set of positive integers containing $1$. If the integer set $\pi_e(G)$ is partitioned into $1$, the set of prime numbers $\pi_e'(G)$ and the set of composite numbers $\pi''_e(G)$, then the following result holds:

\begin{thm}[See \cite{45}]\label{2.7}
   Suppose that $G$ is a finite group, then $|\pi_e'(G)|\leqslant |\pi_e''(G)|+3$. Moreover if $|\pi_e'(G)|=|\pi_e''(G)|+3$, then $G$ is simple and $G$ can be characterized by its spectrum. Furthermore, $G$ is one of the following simple groups.
   \begin{enumerate}[label={\rm{(\Roman*)}}]
        \item $A_5$, $L_2(11)$, $L_2(13)$, $L_2(16)$, $L_3(4)$ and $J_1$.
        \item $Sz(q)$, where $q=2^{2n+1}$ satisfies that each of $q-1$, $q-\sqrt{2q}+1$, and $q+\sqrt{2q} + 1$ is either a prime number or a product of two distinct prime numbers.
        \item $L_2(2^n)$, where $n$ \textup{(}$n\geqslant 5\textup{)}$ is an odd prime and satisfies both $\frac{2^n+1}{3}$ is prime and $2^n-1$ is either a prime number or a product of two distinct prime numbers.
        \item $L_2(3^n)$ where $n$ is an odd prime and satisfies both $\frac{3^n+1}{4}$ is a prime and $\frac{3^n-1}{2}$ is either prime or a product of two distinct prime numbers.
        \item $L_2(5^n)$, where $n$ is an odd prime satisfying both $\frac{5^n-1}{4}$ and $\frac{5^n+1}{6}$ are primes.
        \item $L_2(p)$, where $p$ is a prime greater than $13$ and one of the following holds:
        \begin{enumerate}[label={\rm{(\roman*)}}]
           \item $\frac{p-1}{4}$ and $\frac{p+1}{6}$ are primes.
           \item $\frac{p-1}{6}$ and $\frac{p+1}{4}$ are primes.
        \end{enumerate}
    \end{enumerate}
\end{thm}

\begin{pro}\label{2.8}
In the cases \textup{(II)-(VI)} of the aforementioned theorem, is the number of simple groups satisfying the respective conditions finite or infinite?
\end{pro}

Let $\iota(G)$ denote the set of $pq$ ($p\neq q$) type numbers in the set $\pi_e(G)$. Then the following result holds:

\begin{thm}[See \cite{44}]\label{2.9}
Let $G$ be a finite group, then $|\pi(G)|\leqslant |\iota(G)|+4$. If $|\pi(G)|= |\iota(G)|+4$ then $G$ is simple and $G$ can be characterized by its spectrum.
\end{thm}

From Theorem \ref{2.7} (Theorem \ref{2.9}) we have the number $|\pi(G)|$ of prime divisors of $|G|$ is constrained by $|\pi^{''}_e(G)|$ ($|\iota(G)|$).

If the ``order of a group" is square-free, then the group is solvable (metacyclic group). With respect to this conclusion and the property of $A_5$ (Theorem \ref{1.3}), the author posed the following question:

\begin{pro}\label{2.10}
    Study on finite groups whose ``element orders are square-free".
\end{pro}

Without attaching any additional conditions, what other properties can the spectrum of a finite group possess? That is the Problem \ref{2.2}: What kind of sets of numbers can become a set of element orders $\pi_e(G)$?

After the publication of the literature \cite{160}, Brandl wrote to the author and posed two questions, one of which is: If the orders of elements in a finite group $G$ are consecutive integers, that is, $\pi_e(G)= \{1, 2, 3,\cdots, n\}$, then what is the largest possible value of $n$ that can continue consecutively? Obviously, this question imposes a ``consecutive " condition on the spectrum of the group. In order to solve this question, the author considered prime graph components of finite groups and found the key reference \cite{201} to resolve this problem (Refer to \cite{83, 102} for the relevant literature and to \cite{103} for the most comprehensive description of the prime graph components).

\begin{defi}[See \cite{201}]\label{2.11}

Let $G$ be a finite group. The simple graph $\Gamma(G)$ of $G$ is given by $\pi_e(G)$ as follows: This graph has vertex set $\pi(G)$, and two vertices $p$ and $q$ are adjacent if and only if $pq\in \pi_e(G)$, denoted by $p\sim q$. Denote the connected components of the graph $\Gamma(G)$ by $t(G)$, and using $\pi_i=\pi_i(G)$ $(i=1,2,\cdots, t(G))$, we denote the $i$th connected component of $\Gamma(G)$. If the order of $G$ is even, denote the component containing $2$ by $\pi_1(G)$. This graph $\Gamma(G)$ is usually called the prime graph of group $G$. It was proposed by K. W. Gruenberg and O. H. Kegel in 1975, hence it is also said to be the Gruenberg-Kegel graph of group $G$ \textup{(}abbreviated as prime graph \rm{GK(G)}\textup{)}.
\end{defi}

It is evident from the definition of the prime graph that $\pi_e(G)$ determines $\Gamma(G)$. However, conversely, the same $\Gamma(G)$ can correspond to different sets $\pi_e(G)$. Therefore, using $\Gamma(G)$ to study finite groups, especially finite simple groups, has become a new research topic.

\begin{thm}[See \cite{20}]\label{2.12}
Let $G$ be a finite group. If the elements of spectrum of $G$ are consecutive integers, i.e., $\pi_e(G)=\{1,2,3,\cdots,n\}$ \textup{(}such groups are called finite $\OC_n$-groups\textup{)}, then $n\leqslant 8$. Moreover, these groups have been classified as follows:
\begin{enumerate}[label={\rm{(\Roman*)}}]
\item $n\leqslant 2$ and $G$ is an elementary abelian group.
\item $n=3$ and $G=[N]Q$ is a Frobenius group where either $N\cong Z^t_3$, $Q\cong
      Z_2$ or $N\cong Z^{2t}_2$, $Q\cong Z_3$.
\item $n=4$ and $G=[N]Q$ and one of the following holds:
    \begin{enumerate}[label={\rm{(\roman*)}}]
    \item $N$ has exponent $4$ and class $\leqslant 2$ and $Q\cong Z_3$.
    \item $N\cong Z^{2t}_2$ and $Q\cong \Sigma_3$.
    \item $N\cong Z^{2t}_3$ and $Q\cong Z_4$ or $Q_8$ and $G$ is a Frobenius
          group.
    \end{enumerate}
\item $n=5$ and $G\cong A_6$ or $G=[N]Q$, where $Q\cong A_5$ and $N$ is an
      elementary abelian $2$-group and a direct sum of natural $\SL(2,4)$-modules.
\item $n=6$ and $G$ is one of the following types:
      \begin{enumerate}[label={\rm{(\roman*)}}]
      \item $G=[P_5]Q$ is a Frobenius group, where $Q\cong[Z_3]Z_4$ or $Q\cong \SL(2,3)$ and $P_5\cong Z^{2t}_5$.
      \item $G/O_2(G)\cong A_5$ and $O_2(G)$ is elementary abelian and a direct sum of natural and orthogonal $\SL(2,4)$-modules.
      \item $G=\Sigma_5$ or $G\cong \Sigma_6$.
      \item $n=7$ and $G\cong A_7$.
      \item $n=8$ and $G=[\PSL(3,4)]\langle \beta\rangle$, $\beta$ is an unitary automorphism of $\PSL(3,4)$.
      \end{enumerate}
\end{enumerate}
where $[A]B$ denotes the split extension of its normal subgroup $A$ by a complement $B$.
\end{thm}

\begin{cor}\label{2.13}
    Let $G$ be a finite group, then $\pi_e(G)=\{1,2,3,\cdots,7\}$ if and only if $G\cong A_7$.
\end{cor}

For the infinite $\OC_7$ groups $G$, A. S. Mamontov and E. Jabara \cite{120} proved that it is locally finite, that is, any finite set of elements in $G$ generates a finite subgroup, thereby deducing that $G \cong A_7$ (See \cite[Problem A:19.80]{100}).

It is easy to see that $Z^{\infty}_2$ and $[Z^{\infty}_3]Z_2$ correspond to the infinite $\OC_2$ groups and $\OC_3$ groups, respectively. For $4\leqslant n\leqslant 6$, there also exist infinite $\OC_n$ groups.

\begin{pro}\label{2.14}
     What is the maximum value of $n$ for the infinite $\OC_n$ group? The infinite $\OC_n$ group is periodic, meaning that the orders of elements in the group are finite. Is it locally finite? This is a special case for Burnside's problem \textup{(See \cite{72})}.
\end{pro}

The literature \cite{108, 202} respectively discussed finite groups with orders of elements as consecutive odd integers, and finite groups with orders of elements as consecutive integers except for some primes.

\begin{pro}\label{2.15}
  By analogizing Theorem \ref{2.12}, we investigate finite groups whose spectrum consist of segmented consecutive integers. For example, $\pi_e(A_5)=\{1, 2, 3, 5\}$ \textup{(See \cite{184})} and $\pi_e(\PSL_2(7))=\{1, 2, 3, 4, 7\}$ \textup{(See \cite{160})} both represent finite simple groups with spectrum consisting of two segments of consecutive integers. Furthermore, for finite groups with segmented consecutiveness, determine their maximum length of consecutive integers. For instance, in the case of consecutive, which is an $\OC_n$ group, with a length of $8$.
\end{pro}

Notice that the order of an element in a group is the order of cyclic subgroup. Yanquan Feng \cite{51} and the author \cite{171} provided finite groups with orders of abelian subgroups and orders of proper subgroups are consecutive integers, respectively. Furthermore, since the order of an element in a group is a conjugation invariant. Guohua Qian \cite{150} investigated the finite groups with consecutive nonlinear character degrees.

\section{Characterizing finite simple groups using the order of group and the set of element orders}

According to the classification theorem for finite simple groups (abbreviated CFSG), every finite simple group is isomorphic to one of the following (See \cite{54}):
\begin{enumerate}[label={(\roman*)}]
    \item the groups of prime order;
    \item the alternating groups of degree at least $5$;
    \item the simple classical groups;
    \item the simple exceptional groups of Lie type;
    \item the 26 sporadic groups.
\end{enumerate}

From 1987 to 2003, the author and his coauthors successively demonstrated Conjecture \ref{1.1} to be true for all finite simple groups except the families $B_n(q)$, $C_n(q)$, and $D_n(q)$ (where $n$ is even) (See \cite{32,165,172,180,181,182,205}).  In 2009, this conjecture was also proved for the families $B_n(q)$, $C_n(q)$, and $D_n(q)$ (where $n$ is even) (See \cite{195}). As a result, the aforementioned conjecture has been proved and established a theorem, which states that all finite simple groups can be determined by the ``order of group" and the ``set of element orders" (Simplified as the ``two orders"). This theorem is very helpful in the field of
computational complexity theory (See \cite{BJSP, GJALM}).

For the $26$ sporadic simple groups, their ``two orders" are very clear, which can be studied only by using the prime graph components of finite groups in \cite{201}.
In fact, these $26$ sporadic simple groups unless $J_2$ can all be characterized only by the ``set of element orders" (See Theorem \ref{2.5}(2)).

For alternating groups with degree greater than or equal to $5$, a series of lemmas are deduced in \cite{182} from the order of the group $A_n$, which is $\frac{n!}{2}$.
In particular, the Stirling formula from number theory is used to prove the conclusion without using the prime graph components of finite groups in \cite{201}. In fact, all alternating groups of degree at least $5$, except $A_6$ and $A_{10}$, can be characterized only by using the ``set of element orders" (See Theorem \ref{2.5}(1)).

For the simple classical groups and exceptional simple groups of Lie type, the author first used the classification theorem for simple groups to solve a problem which was posed by E. Artin in 1955: Determining all finite simple groups with the order of Sylow subgroups greater than $|G|^{1/3}$ (See \cite{14}), i.e., the following two lemmas are proved (see \cite{168}):

\begin{lemma}\label{3.1}
Let $G$ be a finite simple group and if $|G|=p^km$, where $p$ does not divide $m$, $p$ is an odd prime and $|G|<p^{3k}$, then $G$ is one of the following groups:
\begin{itemize}
    \item [\rm(1)] A simple group of Lie type in characteristic $p$;
    \item [\rm(2)] $A_5$, $A_6$ and $A_9$;
    \item [\rm(3)] $L_2(p-1)$ where $p$ is a Fermat prime, $L_2(8)$ and $U_5(2)$.
\end{itemize}
\end{lemma}

\begin{lemma}\label{3.2}
    Let $G$ be a finite simple group and if $|G|=2^km$ with $m$ odd, and $|G|<2^{3k}$, then $G$ is one of the following:
    \begin{itemize}
        \item [\rm{(1)}] A simple group of Lie type in characteristic $2$;
        \item [\rm{(2)}] $L_2(r)$ where $r$ is Fermat prime or Mersenne prime;
        \item [\rm{(3)}] $A_6$, $U_3(3)$, $A_9$, $M_{12}$, $U_3(4)$, $A_{10}$, $M_{22}$, $J_2$, $HS$,
              $M_{24}$, $Suz$, $Ru$, $Fi_{22}$, $Co_2$, $Co_1$ and $B$.
    \end{itemize}
\end{lemma}

Above lemmas and the order of group narrow down the scope of the research to simple group of Lie type. Furthermore, upon deeper analysis, Conjecture \ref{1.1} has been proved to be true for all cases except groups $B_n(q)$, $C_n(q)$, and $D_n(q)$ (where $n$ is even).

In fact, for exceptional simple groups of Lie type, Suzuki-Ree groups ${}^2B_2(2^{2n+1})$ ($n\geqslant 1$) (See \cite{169}), ${}^2G_2(3^{2n+1})$ ($n\geqslant 1$) (See \cite{21}), ${}^2F_4(2^{2n+1})$ ($n\geqslant 1$) (See \cite{46}), and $G_2(q)$ (See \cite{196}), $E_8(q)$ (See \cite{104}) and $F_4(2^m)$ (See \cite{31}), and others can be characterized by only using the ``set of element orders" (See \cite[Table 8]{67}).

For classical groups $B_n(q)$ and $C_n(q)$, which have the same order, does there exist some $n$ and prime power $q$ such that their spectra are also the same? Thereby constructing a counterexample to Conjecture \ref{1.1}. Shi \cite{176} and Grechkoseeva \cite{65} almost simultaneously considered this question and proved, using different methods, that for all $n$ and prime power $q$, $\pi_e(B_n(q))\neq \pi_e(C_n(q))$, thus such a counterexample does not exist.

To finally prove Conjecture \ref{1.1}, restrictions on nonabelian composition factors of finite groups have the same spectrum as $B_n(q)$, $C_n(q)$, and $D_n(q)$ (where $n$ is even) in \cite{194}. Finally, Conjecture \ref{1.1} was proved in \cite{195} by using the three theorems of \cite{194}. Consequently, Conjecture \ref{1.1}, which characterizes finite simple groups by using ``two orders", has been proved and established as a theorem.

Numbers and sets of numbers play an important role in mathematics. The concept of ``two orders" was introduced first by the author and has now been extensively studied in the field of group theory. Simple groups are complex, the concept of ``two orders" is the simplest concept. The validation of Conjecture \ref{1.1} establishes a connection between the terms ``simple" and ``complex". Of course, its proof used the classification theorem for finite simple groups.

Except for all finite simple groups we can characterize the following groups using only ``two orders", Jianxing Bi proved the following \cite{17}:

\begin{thm}
    Let $G$ be a finite group and $S_n$ be a symmetry group. Then $G\cong S_n$ if and only if \textup{(1)} $\pi_e(G)=\pi_e(S_n)$; \textup{(2)} $|G|=|S_n|$.
\end{thm}

And Yanwei Gao and Hongping Cao \cite{Gao} proved the following theorem:

\begin{thm}
    Let $G$ be a finite group and $\Aut(G)$ be an automorphism group of sporadic simple group $S$. Then $G\cong \Aut(S)$ if and only if \textup{(1)} $\pi_e(G)=\pi_e(\Aut(S))$; \textup{(2)} $|G|=|\Aut(S)|$.
\end{thm}

In fact, in \cite{Gao} the above condition (1) is weakened to $m_i(G)=m_i(\Aut(S))$ for $i=1,2,3$, where $m_1(G)$, $m_2(G)$ and $m_3(G)$ donote the largest, the second largest and the third largest element order of $G$, respectively.

All finite simple groups can be characterized by ``two orders", but it fails to characterize some groups with small orders. For example: dihedral group $D_8$ with $8$ elements and the quaternion group $Q_8$, obviously, $\pi_e(D_8)=\pi_e(Q_8)=\{1,2,4\}$ but $D_8$ and $Q_8$ are not isomorphism.

\begin{pro}
Except all finite simple groups and some automorphism groups of finite simple groups, what kind of finite groups can be characterized by their two orders?
\end{pro}

Next, we will give another application of Conjecture \ref{1.1} as following:

Let $G$ be a finite group and let $B(G)$ be a Burnside ring of $G$. T. Yoshida posed the following open problem: Let $G$ and $S$ be two finite groups, does $B(G)\cong B(H)$ implies $G\cong H$ (See \cite [Page 340, Proble 2]{211})? It was proved in \cite[Corollary 5.2]{101} that the spectrum of a finite group is determined by its Burnside ring, so we have the following application:

\begin{cor}\label{3.3}
    Let $G$ be a finite group, $S$ a finite simple group. If $B(G)\cong B(S)$ then $G\cong S$.
\end{cor}

Thus, simple groups are recognizable by Burnside ring. In other words, simple groups can be characterized using their Burnside ring.

\begin{pro}\label{3.4}
Is it possible to prove Conjecture \ref{1.1} without using the classification theorem for finite simple groups?
\end{pro}

For a few simple groups with small orders, it is possible to provide a proof without using the classification theorem for simple groups, for example $A_5$ (See \cite{163}). However, it seems impossible to provide a proof without using the classification theorem for all finite simple groups.

\begin{pro}\label{3.5}
Weakening the condition of ``two orders", is it possible to provide a characterization for all finite simple groups? In other words, can they be characterized by the ``group order" and the ``set of orders of certain elements", or by the ``orders of Hall subgroups in a group" and the ``set of orders of elements" for all simple groups?
\end{pro}

For a prime power $q=p^e$, we write $ch(q)=p$, and for a Lie type simple group $G$ defined over $GF(q)$, the characteristic of $G$ denoted by $ch(G)=p$. W. M. Kantor and \'{A}. Seress proved the following theorem (See \cite[Theorem 1.2]{92}):

\begin{thm}\label{3.6}
  Let $G$ and $H$ be simple groups of Lie type of odd characteristic. If $m_i(G)= m_i(H)$ for $i=1$, $2$, $3$, then characteristics of $G$ and $H$ are same, where  $m_1(H)$, $m_2(H)$ and $m_3(H)$ denote the largest, the second largest and the third largest element order of $G$, respectively.
\end{thm}

This gives rise to the following question: Is it possible to characterize simple group using the ``group order" and the ``set of orders of certain elements"? In fact, the following theorem has been proved in \cite{75,76,107,212}:

\begin{thm}\label{3.7}
  Let $G$ be a group. $S$ is one of the following simple groups:
  \begin{itemize}
      \item [\rm{(1)}] Sparodic simple groups \textup{(See \cite{76})};
      \item [\rm{(2)}] $L_2(p)$ where $p\neq 7$ is prime \textup{(See \cite{107})};
      \item [\rm{(3)}] Simple $K_4$-groups with type $L_2(q)$ where $q$ is a prime power \textup{(See \cite{75})};
      \item [\rm{(4)}] Simple $K_5$-groups of type $L_3(p)$, where $p$ is a prime number and $(3,p-1)=1$ \textup{(See \cite{212})}.
  \end{itemize}
  Then $G\cong S$ if and only if $|G|=|S|$ and $m_1(G)=m_1(S)$, where $m_1(G)$ denotes the largest element order of $G$.
\end{thm}

If we add the second largest element order or the third largest element order of group $G$, then it can also be used to characterize more simple groups (See \cite{74}).

Let $G$ be a finite group, set Van($G$)=$\{g\in G|\text{ there exists } \chi \text{ such that } \chi(g)=0 \}$, where $\chi$ is an irreducible complex character
of $G$, Vo($G$) denotes the set of element orders of Van($G$). Obviously, Vo($G$) is a subset of the set of element orders $\pi_e(G)$ of $G$. The following conjecture was posed in \cite [Problem 19.30]{100}:

\begin{con}\label{3.8}
  Let $G$ be a finite group and $S$ a finite simple group, then $G\cong S$ if and only if \rm{(1)} Vo($G$)=Vo($S$); \rm{(2)} $|G|=|S|$.
\end{con}

The conjecture mentioned above  as a question was also put forward in Jinshan Zhang's PhD thesis (See \cite{217}). We refer to the literature \cite{158} and footnote \footnote{Q F. Yan, Z C. Shen, J S. Zhang, et al. A new characteristic of sporadic simple groups $J_1$ and $J_4$. Submitted} for its recent studies.

In fact, the group order, element order, the number of elements of the same order, the length of conjugacy class, and the degree of character, and others, are conjugation invariants. In addition, order of  normalizer of a Sylow subgroup, order of maximal abelian (solvable) subgroup, the index of maximal subgroup, the number of Sylow subgroups, and so on, are important quantities in the study of groups. Study the structure of groups, especially finite simple groups, starting from these important quantities is a broad and meaningful subject. For example, study the structure of finite groups whose lengths of conjugacy classes are all prime powers.

Qian\cite{151} and Qian, et al.\cite{154} defined codegree of the character and gave a connection between the element order and the degree of character. 

Next, we consider the characterization of simple groups by using ``the order of Hall subgroup of a group" and the ``set of element orders". Firstly, the ``order ($|G|_2$) of Sylow $2$-subgroups" of a group $G$ and the ``set of element orders" were used to characterize alternating and sporadic simple groups. We have the following theorems:

\begin{thm}\label{3.9}
    Let $G$ be a finite group and $S$ an alternating simple group. Then $G\cong S$ if and only if \textup{(1)} $|G|_2=|S|_2$; \textup{(2)} $\pi_e(G)=\pi_e(S)$.
\end{thm}

\begin{thm}\label{3.10}
   Let $G$ be a finite group and $S$ a sporadic simple group. Then $G\cong S$ if and only if \textup{(1)} $|G|_2=|S|_2$; \textup{(2)} $\pi_e(G)=\pi_e(S)$.
\end{thm}

\begin{lemma}[See \cite{56}]\label{3.11}
Suppose $G$ is a finite group with $\pi_e(G)=\pi_e(S)$, where $S$ is an alternating simple group $A_n$ and $n\geqslant 5$, $n\neq 6$, $10$. Then $G\cong S$.
\end{lemma}

The authors of \cite{123, 148, 163} earlier studied the characterization of alternating groups and symmetric groups using the set of element orders. The subsequent studies in \cite{56,106,213} the authors further investigated the characterization of alternating groups by using the set of element orders. It was pointed out in \cite{187} that $A_6$ and $A_{10}$ are unrecognizable.

\begin{lemma}[\cite{20}]\label{3.12}
  Let $G$ be a finite group with $\pi_e(G)=\pi_e(A_6)=\{1,2,3,4,5\}$, then $G\cong A_6$ or $G=[N_1]Q$ where $Q\cong A_5$, $N_1$ is an elementary abelian $2$-group and a direct sum of natural $\SL(2,4)$-modules.
\end{lemma}

\begin{lemma}[See \cite{187}]\label{3.13}
    Let $G$ be a finite group with $\pi_e(G)=\pi_e(A_{10})=\{1,2,\cdots,$ $10,12,15,21\}$, then $G\cong A_{10}$, or $G=[A]C$ where $A$ is an Abelian $\{3,7\}$-group, $C=C_G(t)=[\langle t\rangle ]S_5$, where $t$ is an involution and $a^t=a^{-1}$, $a\in A$. And the Sylow $2$-subgroup of $C$ is a generalised quaternion group $Q_{16}$ with $16$ elements.
\end{lemma}

\begin{lemma}[See \cite{134,170}]\label{3.14}
   Let $G$ be a finite group with $\pi_e(G)=\pi_e(M)$ where $M$ is a sporadic simple group distinct from $J_2$, then $G\cong M$.
\end{lemma}

It was proved in \cite{170} that except for $Co_2$ and $J_2$, the remaining $24$ sporadic simple groups can all be characterized by the set of element orders. Mazurov and Shi \cite{134}  proved that $Co_2$ can be characterized by using the ``set of element orders", however if $\pi_e(G)=\pi_e(J_2)=\{1,2,\cdots 8,10,12,15\}$, then $G\cong J_2$, $S_8$, or $G\cong [N]A_8$ is a split extension of a $2$-group $N$ with order $2^{6t}$ ($t=1,2,\cdots $) by $A_8$.
Obviously, from the above arguments and Lemma \ref{3.11}, we can deduce that Theorems \ref{3.9} and \ref{3.10} hold.

\begin{pro}\label{3.15}
Let $G$ be a finite group and $S$ a finite simple group. If \textup{(1)} $|G|_{\pi}=|S|_{\pi}$ where $|G|_{\pi}$ denotes the order of $\pi$-Hall subgroup of $G$, $\pi\neq \pi(G)$; \textup{(2)} $\pi_e(G)=\pi_e(S)$, can we prove that $G\cong S$?
\end{pro}

In 2007, Mazurov posed the following conjecture at International Algebraic Conference held in St. Petersburg: Let $G$ be a finite group, and let $L$ be an alternating group of sufficiently large degree or a simple group of Lie type of sufficiently large Lie rank. If $G$ is isospectral to $L$, then $G$ is an almost simple group with socle isomorphic to $L$. In fact, if $L$ is a finite simple group and $G$ is isospectral to $L$, then either $G$ is solvable, in which case $L$ is one of the following cases: $L_3(3)$, $U_3(3)$, and $S_4(3)$; or $G$  has exactly one nonabelian composition factor (See \cite{228}, \cite[Corollary 7.3]{197} and \cite[Theorem 2]{56}). The conjecture proposed by Mazurov was finally proved to be valid in \cite{68}. By using these results and comparing $|G|_{\pi}$ and $|M|_{\pi}$, Problem \ref{3.15} can be answered.

\section{Thompson Conjecture and Thompson Problem}

\subsection{The origin and early work of the Thompson Conjecture}\

As mentioned earlier, the author posed Conjecture \ref{1.1} to Thompson in 1987, and Thompson proposed the following conjecture in his response letter in 1988
(Thompson Conjecture (1988)): Let $G$ and $M$ be finite groups, set $N(G)=\{n\in\mathbb{Z}^+ \mid G \text{ has conjugacy class } C, |C| = n \}$. Suppose $M$ is a non-abelian simple group, and the center of $G$ is trivial. If $N(G)=N(M)$, then $G$ and $M$ are isomorphic.

Guiyun Chen \cite{33}  proved that if $M$ is a sporadic simpe group then Thompson conjecture is correct. Subsequently, he finished his PhD thesis titled ``On Thompson Conjecture". The author (See \cite{34} ) based on the foundation of prime graph components in \cite{201}, a definition of the order components was provided as follows:

\begin{defi}\label{4.1}
Let $G$ be a finite group, $\pi_1$, $\pi_2,\cdots, \pi_t$ are prime graph components of $G$, where $t=t(G)$ is the number of prime graph components of $\Gamma(G)$. Assume $|G|=m_1m_2\cdots m_t$ where $\pi_i=\pi(m_i)$. Then $m_1$, $m_2, \cdots, m_t$ are called order components of $G$. Set \textup{OC(G)}$=\{m_1, m_2,\cdots, m_{t(G)}\}$, the set of order components of $G$.
\end{defi}

Thompson Conjecture was proved to be true for all sporadic simple groups (See \cite{33,34}). Its validity was also proved for finite group having at least three prime graph components (See \cite{35,36}). For the case of two prime graph components, we can use the same method to prove it. For instance, it was proved that Thompson Conjecture is correct for the group ${}^3D_4(q)$ (See \cite{37}).

\subsection{The case of groups with connected prime graph}\

In 2009, Vasil'ev proved that Thompson Conjecture is valid for the simple groups $A_{10}$ and $L_4(4)$ with connected prime graph(See \cite{193}). Subsequently, N. Ahanjideh proved that Thompson Conjecture is also correct for simple groups of Lie type $A_n(q)$ (See \cite{5}), $B_n(q)$ (See \cite{8}), $C_n(q)$ (See \cite{6}), $D_n(q)$ where $n\neq 4,8$ (See \cite{7}), ${}^2A_n(q)$ (See \cite{9}) and ${}^2D_n(q)$ (See \cite{10}) in a series of articles. For other exceptional simple groups of Lie type, Thompson Conjecture was established the validity in the literature \cite{62,84,97,98,206}. The authors \cite{13,55,57,58,59,60,110,204}  proved that Thompson Conjecture is true for alternating simple groups. Finally, I. B. Gorshkov \cite{61} proved the validity of Thompson Conjecture for the remaining cases $D_4(q)$ and $D_8(q)$. Therefore, Thompson Conjecture was completely proved and was became a theorem.

\subsection{Thompson Problem}\

Thompson problem posed by Thompson in a letter to the author in 1987. Let $G_1$ and $G_2$ be groups of the same order type (as defined in Definition \ref{1.2}). Suppose $G_1$ is solvable, is it true that $G_2$ is also necessarily solvable? For groups $G$ of even order, we can not determine the solvability of $G$ by using the order of $G$ (See \cite{50}). But we
may judge it using the same order type of $G$ if the answer of Thompson Conjecture is positive. In particular, Thompson wrote in the letter: I have talked with several mathematicians concerning groups of the same order type. The problem arose initially in the study of algebraic number fields, and is of considerable interest.
At the same letter, Thompson provided the following example of nonsolvable groups with the same order type in the letter:
$$G_1=2^4:A_7, \hspace{.5cm} G_2=L_3(4):2_2.$$
Both of them are maximal subgroups of $M_{23}$, where ``$:$" denotes semidirect product, we refer to page 72 and page 23 of \cite{41} for the maximal subgroups of $M_{23}$ and the automorphism group of $L_3(4)$, respectively.

The author gave a lecture titled ``Thompson Problem and Thompson Conjecture" at the international conference ``Algebra and Mathematical Logic" held in Russia in September 2011, commemorating the 100th anniversary of V. V. Morozov's birthday. During the conference, the author presented an overview of the aforementioned problem and conjecture.

\subsection{Order equation}\

Let $M_G(n)$ be the set of elements of $G$ with order $n$, the elements of $G$ can be partitioned into the elements of same order, that is,
 $G=\bigcup M_G(n)$ where $n\in \pi_e(G)$.
Since every element of order $n$ must be contained in some cyclic subgroup of order $n$, we have $|M_G(n)|=V_n(G)\phi(n)$, where $V_n(G)$ is the number of cyclic subgroups of order $n$, $\phi(n)$ is Euler totient function.
Hence we can obtain the following equation, $|G| = \Sigma V_n(G)\phi(n)$, $n\in \pi_e(G)$.
Here we denote the ``order equation" of $G$ as $Ord(G)$. That is, if two order equations of finite groups $G$ and $H$ are same,
then $\pi_e(G)=\pi_e(H)$ and $V_d(G)=V_d(F)$ for any $d\in \pi_e(G)$. It is easy to see that same order equation is equivalent to the same order type of Thompson Problem.

Let $G$ and $H$ are two groups of same order type, it is not hard to prove that if $H$ is nilpotent then $G$ is also nilpotent and if $H$ is supersolvable then $G$ is solvable (See \cite{203}). Of course, if $H$ is a finite nonabelian simple group, then $G\cong H$ is nonsolvable because Conjecture \ref{1.1} has been proved. It was proved that if $H$ is solvable and the prime graph of $H$ is disconnected, then $G$ is solvable (See \cite{156,157}). The following results were proved in papers \cite{156,157}:

\begin{lemma}\label{4.2}
Let $G$ be a finite group and $H$ a Frobenius \textup{(}resp. $2$-Frobenius\textup{)} group. If
$G$ and $H$ have same order equation, then $G$ is a Frobenius \textup{(}resp. $2$-Frobenius\textup{)} group. Moreover, if $H$ is solvable, then $G$ is also solvable.
\end{lemma}

From Lemma \ref{4.2} we deduce the following theorem:

\begin{thm}\label{4.3}
    Let $G$ and $H$ be any two finite groups of the same order type. If $H$ is solvable and its prime graph is disconnected, then $G$ is solvable.
\end{thm}

 Note that orders of two elements in a group are equal if they are conjugate. Consequently, ``same order class" is the union of several ``conjugacy classes". Conversely, two elements of the ``same order" might not be ``conjugate". The series of books \textit{Unsolved Problems in Group Theory} records the following question: Suppose that, in a finite group $G$, each two elements of the same order are conjugate. Is then $|G|\leqslant 6$ (see \cite[Problem 7.48]{100})? It was proved in \cite{49,52,191,215} successively using the classification for finite simple groups (CFSG).

\begin{pro}\label{4.4}
  Without using classification for finite simple groups, can we prove that if each two elements of the same order are conjugate in finite group $G$, then $|G|\leqslant 6$?
\end{pro}

Thompson Problem investigates the structure of groups from the ``same order" perspective, while Thompson Conjecture studies groups from the ``conjugacy" perspective. These two have both distinctions and connections. It is well-known that conjugacy classes play a significant role in the study of finite group structure. It was introduced in \cite{30} that the influence of the conjugacy class size or the number of conjugacy classes on finite group structure. $A_5$ was characterized by using the number of elements of the same order (See \cite{155}), that is, the following theorem was proved:

\begin{thm}\label{4.5}
A group $G$ is isomorphic to $A_5$ if and only if $G$ has the same order type as $A_5$ \textup{(} the order type of $A_5$ is \{1, 15, 20, 24 \} \textup{)}.
\end{thm}

Now denote $nse(G)=\{m_k|k\in \pi_e(G)\}$, where $m_k$ denotes the number of elements of order $k$ in $G$, then the above theorem can be stated as follows: A group $G$ is isomorphic to $A_5$ if and only if $nse(G)=nse(A_5)$. For some alternating simple groups and linear simple groups, there is a lot of literature about characterizing these simple groups using $nse(G)$, see \cite{11,15}.

The following result was obtained in \cite{155}:

\begin{prop}\label{4.6}
    Let $G$ be a group \textup{(}not assume finite\textup{)}  containing more than two elements. If the number of elements with the same order in $G$ is finite, and the maximal number is $s$, then $|G|\leqslant s(s^2-1)$.
\end{prop}

Thompson Problem investigates that using ``same order" to judge the solvability of groups. The condition ``same order" cannot be changed to ``conjugacy" here. In fact, we have the following conclusion:

\begin{thm}[See \cite{145}]\label{4.7}
There exist two finite groups $G$ and $H$ such that $N(G)=N(H)$, i.e., the set  of their conjugacy class sizes are same, but $G$ is solvable and $H$ is nonsolvable.
\end{thm}

A related work to Thompson Problem is the study of the solvability of groups using the number of elements of maximal order or the type of this number (See \cite{38,48,77,87,88,89,90,91,189,207}).

Next, we continue discussing Thompson problem and give the following definitions.

For each finite group $G$ and each integer $d\geqslant 1$, let $G(d)=\{x\in G|x^d=1\}$.

\begin{defi}
$G_1$ and $G_2$ are of the same two orders type if and only if $|G_1|=|G_2|$ and $\pi_e(G_1)=\pi_e(G_2)$.
\end{defi}

\begin{defi}
$G_1$ and $G_2$ are of the same Burnside ring type if and only if $B(G_1)\cong B(G_2)$.
\end{defi}

The authors \cite{LIY} proved that if $G_1$ and $G_2$ are of the same Burnside ring type, then $G_1$ and $G_2$ are of the same order type (see Def. 1.3). Thus $G_1$ and $G_2$ are of the same two orders type immediately. But the reverse side is not correct. We have the following examples:

\begin{exa}
    Let $G_1=L_3(4):2_2$ and $G_2=2^4:A_7$. Then $G_1$ and $G_2$ are of the same order type, but $B(G_1)\not\cong B(G_2)$.
\end{exa}

\begin{exa}\label{4.12}
    Let $G_1=A_5\times A_5\times A_5$, $G_1$ is not solvable and $|G_1|=60^3$ and $\pi_e(G_1)=\{1, 2, 3, 5, 6, 10, 15, 30\}$ since $\pi_e(A_5)=\{1, 2, 3, 5\}$. Let $G_2=C_{30}\times C_{30}\times C_{30}\times C_2\times C_2\times C_2$.
    Then $|G_2|=60^3$ and $\pi_e(G_2)=\{1, 2, 3, 5, 6, 10, 15, 30\}$. $G_1$ and $G_2$ are of the same two orders type, but $G_2$ is solvable. Moreover, $G_1$ and $G_2$ are of not the same order type.
\end{exa}

Noting Thompson problem, which states that if $G_1$ and $G_2$ are of the same order type, and $G_1$ is a solvable group, is $G_2$ also solvable? If we strengthen ``same order type" to ``same Burnside ring type", then the answer is positive (See \cite{DA}). However, if we weaken ``same order type" to ``same two orders type", then the above Example \ref{4.12} is a  counterexample.

Recently, the author of \cite{229} construct two finite groups of size $2^{365}· 3^{105}· 7^{104}$: a solvable group $G$ and a non-solvable group $H$, such that for every integer $n$ the groups have the same number of elements of order $n$. It gives a negative answer to Thompson problem. As a  much smaller counterexample to Thompson problem, reference \cite{231} constructed two same order type groups of order $2^{13}· 3^4· 7^3$, one solvable and the other non-solvable.

\section{Problems related to quantitative characterization}

\subsection{Width of order and width of spectrum in finite groups}\

The fundamental theorem of arithmetic states that every integer $n$ can be represented as $n=p^aq^b\cdots r^c$, where $p$, $q,\cdots, r$ are distinct prime numbers.

The prime factorization of numbers is one of the most important results in mathematics. In the following section, we will study the influence of the number of prime divisors on the structure of groups.

Note that both the order $|G|$ of a finite group $G$ and the order $|g|$ of an element $g$ are numbers in $G$, let $|\pi(G)|$ and $|\pi(g)|$ denote the number of prime factors of $|G|$ and $|g|$, respectively. These are respectively called the width $\omega_0(G)$ of the order of finite group $G$ and the width of an element $g$.

\begin{defi}
Let $G$ be a finite group, we define $\omega_0(G)=|\pi(G)|$ the width of order of $G$ and $\omega_s(G)=max\{|\pi(g)|\big|g\in \pi_e(G)\}$ the width of spectrum of $G$.
\end{defi}

For other sets of conjugation invariants of group $G$, the corresponding width can also be defined.

In August 2020, at the online conference ``Ural Workshop on Group Theory and Combinatorics" held in Russia, the author gave the first presentation and introduced the aforementioned two definitions of width, and discussed the structure of finite groups with a small width, we refer to \cite{179} for some results about the ``width of order". This paper points out that groups with $\omega_0(G)=1$, finite $p$-groups, although the orders of these groups have only one prime factor, their structures are quite complex. Recent Chinese literature can be found in \cite{221,222}. For the especially ``simple" condition, $\pi_e(G)=\{1,3\}$, but its structure is quite complex, and currently, there is no a detailed classification.

When $\omega_0(G)=2$, that is, $|G|=p^aq^b$, it is a special class of solvable groups. In 1904, Burnside \cite{25} provided the first proof of its solvability using representation theory. Later, D. M. Goldschmidt \cite{53} provided a purely group theoretical proof for this theorem when both $p$ and $q$ are odd primes. Then, H. Bender \cite{16} and H. Matsuyama \cite{121}, in 1972 and 1973 respectively, proved the remaining cases, thereby providing an abstract group proof for the $p^aq^b$ theorem.

If $\omega_0(G)=3$, then $G$ can be nonsolvable group, and its chief factors are finite simple groups with $\omega_0(G)=3$. M. Herzog \cite{80} gave all $8$ finite simple groups (simple $K_3$-groups) with $\omega_0(G)=3$, they are $A_5$, $L_2(7)$, $L_2(8)$, $A_6$, $L_2(17)$, $L_3(3)$, $U_3(3)$, and $U_4(2)$.

Using classification for finite simple groups, simple groups with $\omega_0(G)=4$ (simple $K_4$ groups) were successively given in the papers \cite{82,167,200}
\footnote{Reference \cite{200} contains conclusions without specific proofs, and there are omissions (such as $L_3(8)$ and $U_3(7)$) and repetitions (for instance $O_7(2)$ and $S_6(2)$ are the same group) in some places, furthermore, Suzuki simple groups can be completely determined.}
. The further investigation of the simple group with $\omega_0(G)=4$ (simple $K_4$-group) has been discussed in \cite{24,223}, it is one of the following groups:

\begin{enumerate}[label={(\Roman*)}]
    \item $A_n$ ($n=7,8,9,10$), $M_{11}$, $M_{12}$, $J_2$; $L_2(q)$ ($q=16,25,49,81,97,243,577$), $L_3(q)$ ($q=4,5,7,8,17$), $L_4(3)$; $O_5(q)$ ($q=4,5,7,9$), $O_7(2)$, $O^+_8(2)$, $G_2(3)$; $U_3(q)$ ($q=4,5,7,8,9$), $U_4(3)$, $U_5(2)$; ${}^3D_4(2)$, ${}^2F_4(2)'$, $Sz(8)$, $Sz(32)$;
    \item $L_2(r)$, where $r$ is prime and satisfies the following equation: $r^2-1=2^a3^bu$ for $a\geqslant 1$, $b\geqslant 1$, $u>3$ is prime;
    \item $L_2(2^m)$ and the following equations hold: $2^m-1=u$, $2^m+1=3t$ where $m\geqslant 1$, $u$ and $t$ are prime numbers and $t>3$;
    \item $L_2(3^m)$ and the following equations hold: $3^m+1=4t$, $3^m-1=2u$ where $m\geqslant 1$, $u$ and $t$ are odd primes.
\end{enumerate}

Every group of order $p^aq^b$ is solvable, so the number of simple $K_2$ group is zero. However, the number of simple $K_3$ groups is $8$. The following is a conjecture about the number of simple $K_4$ groups.

\begin{con}\label{5.2}
There are infinitely simple $K_4$ groups.
\end{con}

In particular, the author conjectured that there are infinite prime number $r$ such that $\omega_o(L_2(r))=4$. However, it is difficult to prove this result because it is a special case of the unsolved Dickson conjecture (See \cite{47}).
In the following pairs of numbers
$$(x,3x-2),\quad (x,2x+1),\quad (x,4x+1),\quad (x,6x-1),$$
$$(x,6x+1),\quad (x,2x-1),\quad (x,4x-1),\quad (x,8x-1),$$
if it can be proved that for prime number $x$ such that there are infinitely many prime pairs appearing in any of the aforementioned pairs, then the number of simple $K_4$ groups is infinite. Conversely, if for infinite prime numbers $x$, the second position of the above pairs can only appear finite prime numbers, can a contradiction be deduced from this?

Simple groups with $\omega_0(G)=5$ and $\omega_0(G)=6$ were studied in \cite{86,105}.

The following results hold for the width $\omega_0(G)$ of order of $G$: $\omega_0(G)\leqslant |\pi^{''}_e(G)|+3$, and if $\omega_0(G)=|\pi^{''}_e(G)|+3$ then $G$ is simple where $\pi^{''}_e(G)$ denotes the set of composite numbers in the set $\pi_e(G)$ (See Theorem \ref{2.7}); $\omega_0(G)\leqslant |\iota(G)|+4$, and if $\omega_0(G)=|\iota(G)|+4$ then $G$ is simple, where $\iota(G)$ denotes the set of $pq$ ($p\neq q$) type numbers in $\pi_e(G)$ (See Theorem \ref{2.9}).

Similar to Conjecture \ref{5.2}, there also arise the following number theory problems from the classification theorem for finite simple groups:

\begin{pro}\label{5.3}
  Are there infinitely simple groups whose order being a square number?
\end{pro}

In fact, the arising of the above problem is also originated from the classification of simple groups.
In the historical process of researching the classification of simple groups, the work of studying simple groups whose order precisely contains a prime number to the first power played an important role. Brauer and Duan have made outstanding contributions to this topic. Brauer \cite{BR2} pointed out that among simple groups with orders less than $10^9$, only one simple group, namely the order of the classical group $B_2(7)$, does not precisely contain a prime number to the first power. However, determining which simple groups do not have orders precisely containing a prime number to the first power is a challenging problem, and it seems difficult to determine whether the number of such simple groups is finite or infinite. Zhongmu Chen \cite{Chen} proved the following theorem.

\begin{thm}
The order of a finite simple group $G$ is not a $k$-th power, where $k\geqslant 3$. The necessary and sufficient condition for the order of a group $G$ to be a perfect square is that $G$ is a simple group $B_2(p)$ of Lie type, where $p$ is a prime number satisfying
$$p=1+\binom{2n+1}{2}2+\binom{2n+1}{4}2^2+\cdots+\binom{2n+1}{2n}2^n, n\geqslant 1.$$
The order of $B_2(p)$ is $p^4(p^2-1)(p^4-1)/2$. For example, if $n=1$, then $p=7$, and $|B_2(7)|=2^8\cdot3^2\cdot 5^2\cdot 7^4$.
\end{thm}

However, in the above expression for the summation, the number of $n$ for which the sum is a prime number, remains an unresolved problem.

\subsection{Finite groups with a small spectrum width}\

If we consider a group as a whole, its elements can be seen as local, forming an inseparable connection between the local and the global. Our consideration of the element orders is a way to understand the whole from a local perspective. As previously mentioned, for any finite group $G$ and any element $g$ in $G$, the order $|g|$ of $g$ divides the order $|G|$ of $G$. During the author's graduate studies, the author considered such a question: the study of finite groups in which the orders of elements, excluding the identity element, are all prime numbers. This idea has received full approval from the author's advisor, Zhongmu Chen. Furthermore, this work was extended to study groups in which the orders of elements are prime powers, i.e., groups whose spectrum width is $1$. These are the topics of the author's master's thesis and they were studied in \cite{184,185}. Now we say a group $G$ is called prime power order group (briefly, $\EPPO$-group) if $\omega_s(G)=1$. Moreover, a group $G$ is called prime element group (briefly, $\EPO$-group) if its every nontrivial element has prime order.

After completing the master's thesis, the authors \cite{185} found that G. Higman had previously studied such groups and proved the following theorem (See \cite[Theorem 1]{HG}).

\begin{thm}\label{EPPO}
Let $G$ be a soluble $\EPPO$-group. Let $p$ be the prime such that $G$ has a normal $p$-subgroup greater than $1$, and let $P$ be the greatest normal $p$-subgroup of $G$. Then $G/P$ is either
\textup{(i)} a cyclic group whose order is a power of a prime other than $p$; or \textup{(ii)} a
generalized quaternion group, $p$ being odd; or \textup{(iii)} a group of order $p^aq^b$ with cyclic Sylow subgroups, $q$ being a prime of the form $kp^a+1$. Thus $G$ has order divisible by at most two primes, and $G/P$ is metabelian.
\end{thm}

For solvable $\EPPO$-groups, the authors of \cite{185} provided a more detailed structure of groups than Theorem \ref{EPPO}.

\begin{thm}[See Theorem 2.4 of \cite {185}]
Let $G$ be a solvable $\EPPO$-group such that $|\pi(G)|>1$ and $Q$ be a maximal normal $q$-subgroup of $G$ for some $q\in \pi(G)$. Then
\begin{itemize}
\item [\rm{(1)}] Suppose that $G_2$ is not a generalized quaternion group, where $G_2$ is a Sylow $2$-subgroup of $G$. Then $G/Q$ is a metacyclic group. Let $|G|=p^{\alpha}q^{\beta}$. Then $G/Q$ is of order $p^{\alpha}q^{\gamma}$ with $q^{\gamma}|(p-1)$ and
the chief series of $G$ is as follows:
$$\underbrace{q,\cdots,q}_{{\gamma}};\underbrace{p,\cdots,p}_{{\alpha}};q^{b_1},\cdots,q^{b_k};\quad b|b_i, i=1,\cdots,k,\quad \gamma<b,$$
where $p^{\alpha}|(q^b-1)$. If $Q$ is the Sylow $q$-subgroup of $G$, then the chief series of $G$ is as follows:
$$\underbrace{p,\cdots,p}_{\alpha};\underbrace{q^b,\cdots,q^b}_{k}, \quad \beta=kb,$$
and the length of nilpotency class of $Q$ is bounded by $k$.
\item [\rm{(2)}] If the Sylow $2$-subgroups of $G$ are generalized quaternion groups, then $G$ has the following chief series
$$\underbrace{2,\cdots,2}_{\alpha};q^{b_1},\cdots,q^{b_k},$$
where $b_i>1$, $b|b_i$ for $i=1,\cdots, k$, where $b$ is the exponent of $q(\bmod~2^{\alpha-1})$.
\end{itemize}
\end{thm}

For nonsolvable $\EPPO$-groups, based on the result of Suzuki \cite{188}, Brandl \cite{BR} and the authors of \cite{185} obtained the following theorem (See \cite [Theorem 3]{BR}, the Mathieu group $M_{10}$ is omitted, and \cite [Theorem 3.1]{185}).

\begin{thm}\label{nonsolvable}
Suppose that $G$ is a nonsolvable $\EPPO$-group. Then one of the following holds.
\begin{itemize}
\item [\rm{(1)}] $G$ is isomorphic to $\PSL_2(q)$ for $q=7, 9$ or $\PSL_3(4)$ or $M_{10}$.
\item [\rm{(2)}] There exists a normal $2$-subgroup $T$, $T$ is an elementary abelian group such that $G/T$ is isomorphic to
\begin{enumerate}[label={\rm{(\roman*)}}]
\item $\PSL_2(q)$ for $q=5, 8, 17$. In this case, the class length of $T$ is not greater than $2$ and $|T|-1$ is divisible by $3\cdot5$, $3^2\cdot7$, $3^2\cdot 17$, respectively.
\item $Sz(2^3)$ with $5\cdot 7\cdot 13\big| (|T|-1)$.
\item $Sz(2^5)$ with $5^2\cdot 31\cdot 41\big|(|T|-1)$.
\end{enumerate}
\end{itemize}
\end{thm}

From the above Theorem \ref{nonsolvable} we can get the following result:
\begin{cor}
    Let $G$ be a finite group. Then $G\cong M_{10}$ if and only if $\pi_e(G)=\{1,2,3,4,5,8\}$.
\end{cor}

About Theorem \ref{nonsolvable} the authors of \cite{29} have summerized in their paper all known results about non-solvable EPPO-groups in the form "if and only if" (See \cite[Theorem 1.7]{29}).

We also obtained the following interesting result (See \cite[Theorem 2.1]{185}):
\begin{thm}
    A solvable $\EPPO$-group $G$ is an $M$-group. That is, every irreducible representation of $G$ can be induced from a linear representation of a subgroup of $G$, i.e., a monomial representation.
\end{thm}
Note that the abstract group definition of an $M$-group is an unresolved question.

The following result was proved for $\EPPO$-group (See \cite[Theorem 1.4]{185}):

\begin{property}\label{5.4}
Let $G$ be an $\EPPO$-group and $H$ be a subgroup of $G$ such that $(|H|, d)=1$ for a natural number $d>1$. Then $|H|$ divides the number of elements of order $d$ in $G$.
\end{property}

The following result shows that the aforementioned property is a characteristic property of $\EPPO$-groups.

\begin{thm}[See \cite{26}]\label{5.5}
 Let $G$ be a finite group and let $H<G$ such that $(|H|,d)=1$ where $1\neq d\in \pi_e(G)$. Then $G$ is a $\EPPO$-group if and only if $|H|$ divides the number of elements of order $d$ in $G$.
\end{thm}

The study of finite $\EPPO$-groups was extended to the infinite case in \cite{79,209}.

\begin{pro}\label{5.6}
 Study the structure of finite groups $G$ with $\omega_s(G)=2$ and $\omega_s(G)=3$.
\end{pro}

T. M. Keller and A. Moret\'{o} \cite{96} studied finite groups $G$ with $\omega_s(G)=2$, that is, there exists $pq\in \pi_e(G)$ such that $pqr\not\in \pi_e(G)$ where $p$, $q$, $r$ are different primes. For the special cases, i.e., solvable $\{p,q,r\}$-groups, the following result holds (See \cite[Theorem C]{96}):

\begin{thm}
Let $G$ be a solvable $\{p,q,r\}$-group, $p$, $q$, and $r$ divide $|G|$. If $pqr\not\in \pi_e(G)$ \textup{(}i.e., $\omega_s(G)\leqslant 2$\textup{)}, then $h(G)\leqslant 21$. Furthermore, if $|G|$ is odd, then $h(G)\leqslant 15$, where $h(G)$ is the Fitting height of $G$.
\end{thm}

Theorem $C$ of \cite{96} is a recently obtained result, it can be seen that study the structure of finite groups with $\omega_s(G)=2$ is a rather challenging problem. Qian \cite{152} has recently provided the structure of solvable groups whose spectrum width is $3$. Note that, the classification of finite groups with $\omega_s(G)=1$ ($\EPPO$-group) was accomplished based on the work in \cite{188}. The classification for simple groups might be used to determine the structure of finite groups satisfying $\omega_s(G)=2$ and $\omega_s(G)=3$.

\subsection{Generalizations of related problems}\

Using spectrum to characterize the finite alternating group $A_5$ was a study conducted three to four decades ago. Here are several generalizations from this work.

``Spectrum" is the set of element orders, i.e., it is the set of orders of cyclic subgroups. We can consider the influence of sets of numbers on group structures, as mentioned after Conjecture \ref{3.8}: the set of orders of abelian subgroups, the set of orders of nilpotent subgroups, the set of orders of solvable subgroups (See \cite{19, 43}), the set of orders of the normalizer of Sylow subgroups (See \cite{18} ), the set of indices of maximal subgroups (See \cite{LXH}), etc. The order of element, the number of elements of the same order, the length of conjugacy classe, and the degree of character, and others, are conjugation invariants in $G$. Different sets of numbers can be constructed from these conjugation invariants.

``Characterization" can be ``isomorphic" or ``homomorphic". From the perspective of the ``function $h$", it can be $h(\Gamma)=1$ or $h(\Gamma)=k$ where $k$ is finite, otherwise, it is called unrecognizable.

``Finite simple group" can be a ``simple group" or a nonsolvable group, such as the automorphism groups of some simple groups, it can also be the direct product of simple groups, for example $Sz(2^7)\times Sz(2^7)$ (See \cite{124}) and $J_4\times J_4$ (See \cite{64}).

``Spectrum" represents the research conditions, ``characterization" represents the research methodology, and ``finite simple groups" represent the research subject. Many questions can be raised from various perspectives to conduct research on these three aspects.

\subsection{The quantity relationship between the width of order and the width of spectrum}\

In 1989, the author gave a lecture at The $3^{th}$ National Conference on Algebra in China, and the main content was published in \cite{166}. In this lecture, the author summarized the early research on characterizing simple groups by using the ``order" and introduced Thompson Problem and Thompson Conjecture in \cite{166} (see Sections 4.1 and 4.4 of this article). The author \cite{166} posed the following question (see \cite[Question 4.3] {166}):

Given a positive integer $k$, denote by $|\pi(k)|$ the number of distinct prime divisors of $k$. Let $G$ be a finite group, $n=max\{|\pi(k)|\big | k\in \pi_e(G)\}$, does there exist a function $f$ such that $|\pi(G)|\leqslant f(n)$? In other words, for any finite group, is the width of its order bound by the width of its spectrum?

If $\omega_s(G)=n=1$, then $G$ is a finite $\EPPO$-group, from \cite[Theorem 2.4, 3.1 and 3.2]{185} we deduce that $\omega_0(G)=|\pi(G)|\leqslant 4$. Keller \cite{KTM} proved that if $G$ is solvable and $\omega_s(G)=3$ then $\omega_0(G)\leqslant 8$.

For the above question, Jiping Zhang \cite{216}, C. Bellotti, Keller and T. S. Trudgian \cite{BKT} proved that if $G$ is solvable, then the width $\omega_0(G)$ of the order is bounded by a quadratic function and linear function of the spectrum width $\omega_s(G)$, respectively. That is,

 $$\omega_0(G)\leqslant \frac{\omega_s(G)(\omega_s(G)+3)}{2} \text{ (See \cite{216})},$$
 $$ \omega_0(G)\leqslant 5\omega_s(G) \text{ (See \cite{BKT})}.\qquad \quad \qquad $$

For general group, $\omega_0(G)$ is bounded by a super-exponential function of $\omega_s(G)$. The above result was subsequently improved in the literature \cite{94,95, 144}. Recently, the following theorem was proved in \cite{81}:

\begin{thm}\label{5.7}
Let $G$ be a finite group, then $\omega_0(G)<210\omega_s(G)^4$.
\end{thm}

An analogous result regarding the order of a group and its co-degrees can be found in \cite{210}.

It is easy to get that $\omega_s(G)\leqslant \omega_0(G)$ by the definitions of the order width and the spectrum width. The following result holds for all finite simple groups:

\begin{thm}[See \cite{179}]\label{5.8}
Let $G$ be a finite simple group, then $\omega_s(G)<\omega_0(G)$.
\end{thm}

Naturally, the author of this article put forward a question to investigate the difference between two widths, denoted as $d=\omega_0(G)-\omega_s(G)$. For example, discuss the classification of finite simple groups with $d=1$.

\section{Group, quantity and graph}

\subsection{Group and quantity}\

Groups and quantities are closely connected. The group order can be used to determine some important properties of the group, such as the groups of order $15$ are cyclic, the groups of order $p^2$ are abelian, and the groups of order $135$ are nilpotent, and so on. For the solvability of groups, there are many well-known results were obtained in the research history of group theory, such as Feit-Thompson odd order theorem (See \cite{50}) and Thompson classification theorem for minimal simple group (See \cite{190}). From these, the following theorem can be obtained:

\begin{thm}\label{6.1}
 Let $G$ be a finite group. If $(|G|,2)=1$ or $(|G|,15)=1$, then $G$ is solvable.
\end{thm}


The numbers $2$ and $15$ in the theorem above are called ``solvable coprime number" in lecture \cite{73}, meaning that a group whose order is coprime to $2$ or $15$ is necessarily a solvable group. Moreover, it is obtained in \cite{73} that all solvable coprime numbers must be multiples of either $2$ or $15$. Compared to the study of finite groups using the order of the group, we characterize finite simple groups using the set $\pi_e(G)$ of element orders. In addition, the set $\pi_e(G)$ of element orders can also be used to judge the solvability of $G$.

\begin{defi}\label{6.2}
Let $G$ be a finite group and let $\pi_e(G)$ be the set of element orders in $G$. The set $T$ of numbers is called the set of "intersection empty solvability"
if $\pi_e(G)\cap T=\varnothing$ implies that $G$ is a solvable group.
\end{defi}

\begin{thm}[See \cite{178}]\label{6.3}
Let $G$ be a finite group and let $\pi_e(G)$ be the set of element orders in $G$. If $\pi_e(G)\cap T=\varnothing$ where $T=\{2\}$, $\{3,4\}$ or $\{3,5\}$, then $G$ is solvable. Furthermore, using the set T of "intersection empty solvability" to judge $G$ is solvable or not, only the above three cases \textup{(} $T=\{2\}$, $\{3,4\}$ or $\{3,5\} $\textup{)}.
\end{thm}

\begin{pro}\label{6.4}
Let $G$ be a finite group. Is it possible to use $\pi_e(G)$ to provide sufficient conditions for $G$ to be cyclic, abelian, nilpotent, supersolvable, or an $M$-group, etc.?
\end{pro}

Compared to the quantity problem of the ``order of group", some similar problems can also be posed by using the ``element orders". For example, Lagrange Theorem of finite group $G$ implies that the order of every subgroup of $G$ divides the order of $G$. However, for any divisor of $|G|$, a subgroup of that order may not necessarily exist. In other words, the converse of Lagrange's theorem does not hold.

\begin{defi}\label{6.5}
A finite group $G$ is called $\CLT$ \textup{(}the converse of Lagrange's theorem\textup{)} group if and only if the converse of Lagrange's theorem holds for $G$.
\end{defi}

\begin{defi}\label{6.6}
 A subset $M$ of a set of positive integers is called a closed subset if, for all $k\in M$, all factors of $k$ also belong to $M$.
\end{defi}

Let $G$ be a group and $\pi_e(G)$ be the set of element orders of $G$. It is not all closed subsets of $\pi_e(G)$ can become the set of element orders of a subgroup of $G$.

Compared to $\CLT$ group, we introduce the following definition:

\begin{defi}\label{6.7}
We say that a finite group $G$ is a $\COE$ group if any closed subset $M$ of $\pi_e(G)$, there exists a subgroup $H$ such that $\pi_e(H)=M$.
\end{defi}

\begin{thm}[See \cite{174}]\label{6.8}
Suppose $G$ is a $\COE$ group, then $|\pi(G)|\leqslant 3$, and one of the following holds.
\begin{itemize}
    \item [\rm{(1)}] $G$ is a $p$-group;
    \item [\rm{(2)}] $|G|=p^aq^b$ where $p\neq q$, $ab\neq 0$, and if the exponent of $G$ is $p^mq^n$ then $G$ has an $\EPPO$-subgroup with an exponent of $p^mq^n$.
    \item [\rm{(3)}] $|G|=2^a3^b5^c$ where $abc\neq 0$ and the exponent of $G$ is $30$, $60$ or $120$. Furthermore, $G$ has a proper subgroup that is isomorphic to $A_5$, and $G$ is not an $\EPPO$-group.
\end{itemize}
\end{thm}

\subsection{Group and graph are closely interconnected}\

Gruenberg and Kegel
\footnote{K. W. Gruenberg, O. H. Kegel. Unpublished manuscript, 1975} first introduced the definition of prime graph $\Gamma(G)$ of group $G$ (See Definition \ref{2.11}) and its classification. Subsequently, the literature \cite{83,102,103,201} provided explicit descriptions for the connected components of the prime graphs of the finite simple groups. M. S. Lucido studied the diameter of the prime graph as follows (See \cite{111}):

\begin{defi}\label{6.9}
We say that $$\diam(\Gamma(G))=max\{ d(p,q)|p,q \textit{ in the same connected component of } \Gamma(G)\}$$ is the diameter of the prime graph of a group $G$, where $d(p,q)$ denotes the distance between elements $p$, $q$ if they are in the same connected component of $\Gamma(G)$.
\end{defi}

If $G$ is a finite solvable group then $\diam(\Gamma(G))\leqslant 3$. Lucido \cite{111} proved that for all finite groups $G$, we have $\diam(\Gamma(G))\leqslant 5$ and characterized almost simple groups with $\diam(\Gamma(G))=5$. Recall that a group $B$ is called an almost simple group if there exists a nonabelian simple group $A$ such that $A\leqslant B\leqslant \Aut(A)$.

Finite groups with $\diam(\Gamma(G))=5$ were obtained in \cite{63}.

Note that, when $\diam(\Gamma(G))=1$, the group $G$ is an $\EPPO$-group (See \cite{185}. So, we pose the following question:

\begin{pro}\label{6.10}
Study finite groups with $\diam(\Gamma(G))=2$, $3$, $4$.
\end{pro}

In the case of solvable groups, the literature \cite{96, 152} respectively investigated finite solvable groups with $\diam(\Gamma(G))=2$ and $\diam(\Gamma(G))=3$.

Lucido \cite{112} further discussed the case of a prime graph being a tree based on the work in reference \cite{111}, i.e., a connected graph without loops. The following result was obtained:

\begin{thm}\label{6.11}
If the prime graph of a finite group $G$ is a tree, then $|\pi(G)|\leqslant 8$.
\end{thm}

Lucido \cite{112} also gave an example with $|\pi(G)|= 8$.

Compared to the prime graph $\Gamma(G)$ of $G$, the definition of solvable graph was introduced in \cite{4} as follows:

\begin{defi}\label{6.12}
Let $G$ be a finite group. We define the solvable graph $\Gamma_s(G)$ of $G$ as follows: its vertices are the primes dividing the order $|G|$ of $G$, and vertices $p$ and $q$ are joined by an edge if and only if there exists a solvable subgroup whose order is divisible by $pq$.
\end{defi}

M. Hagie \cite{71} proved that $\diam(\Gamma_s(G))\leqslant 4$.

There is a lot of work in study groups using graphs, P. Erd\"{o}s defined non-commuting graph in 1976 (See \cite{147}).

\begin{defi}\label{6.13}
Let $G$ be a finite group and let $Z(G)$ be the center of $G$. We define the non-commuting graph $\triangledown(G)$ of $G$ as follows: Take $G\backslash Z(G)$ as the vertices of $\triangledown(G)$ and join two distinct vertices $x$ and $y$ if and only if $xy\neq yx$.
\end{defi}

The following problem was conjectured in \cite{1,140}: Let $G$ and $H$ be two nonabelian finite
groups such that $\triangledown(G)=\triangledown(H)$, then $|G|=|H|$. A. R. Moghaddamfar \cite{137} gave some counterexamples to this conjecture. Moreover, the following conjecture was posed in \cite{1}, which is called AAM (Abdollahi-Akbari-Maimani) Conjecture:

\begin{con}\label{6.14}
Let $S$ be a nonabelian simple group and $G$ a group. If $\triangledown(G)=\triangledown(S)$, then $G\cong S$.
\end{con}

It was proved in \cite{199,220} that for alternating simple group $A_{10}$ and $L_4(4)$, whose prime graphs are disconnected, Conjecture \ref{6.14} is valid. In addition, the following result was also proved in \cite{199}:

\begin{prop}\label{6.15}
Let $S$ be a nonabelian simple group and let $G$ be a group with $Z(G)=1$. If $\triangledown(G)=\triangledown(S)$, then $N(G)=N(S)$.
\end{prop}

Because Thompson Conjecture has been proved (See sections 4.1 and 4.2), consequently, AAM Conjecture is correct.

\begin{defi}[OD-characterization of simple group]\label{6.16}
Let $G$ be a finite group and $|G|=p^aq^b\cdots r^c$, where $p<q<\cdots <r$ are primes, and $a$, $b\cdots, c$ are integers. For $s\in \pi(G)$, let $deg(s)=|\{t\in \pi(G)|s\sim t\}|$, which we call the degree of $s$. We also define $D(G)=(deg(p), deg(q),\cdots, deg(r))$, we call $D(G)$ the degree pattern of $G$. A group $M$ is called $k$-fold OD-characterizable if there exist exactly $k$ non-isomorphic groups $G$ such that $|G|=|M|$ and $D(G)=D(M)$. Moreover, a $1$-fold \textup{OD}-characterizable group is simply called an \textup{OD}-characterizable group.
\end{defi}

Extensive research has been conducted by Moghaddamfar, Liangcai Zhang, et al. on OD-characterizable ($k$-fold OD-characterizable) simple group (almost simple group) (See \cite{12,141,142,143,218,219}).

\subsection{Adjacency criterion for the vertices of the prime graph (GK graph) of a finite simple group}\

Both the characterization of finite simple groups by using the set of element orders $\pi_e(G)$ and the aforementioned OD-characterization of finite simple groups are inseparable from the study of the prime graph of groups (GK graph) (See Definition \ref{2.11}). Therefore, studying the adjacency of the vertices in prime graph is a very important subject. The literature \cite{197,198} gave the adjacency criterion for the vertices of prime graphs (GK graph) associated with each nonabelian simple group. Especially in the case of vertex adjacency to vertex ``2" and in the case of ``independent sets" in which vertices are not adjacent to each other. The literature \cite{197,198} provided a large number of graphs and information for all simple groups in tabular form. Reference \cite{198} is a continuation of reference \cite{197}, and it includes corrections to some content in reference \cite{197}.

 \subsection{Unrecognizable finite group by prime graph}\

Finite simple groups $E_6(2)$, $E_6(3)$ and ${}^2E_6(3)$ were characterized by its prime graph $\Gamma(G)$ in \cite{69, 99}. Clearly, there are many simple groups that cannot be only characterized using prime graphs. Thus, P. J. Cameron and N. V. Maslova \cite{29} studied the criteria for finite groups that cannot be characterized using prime graphs and they aimed to study the following questions:

\begin{enumerate}
    \item  Which groups $G$ are uniquely determined by their prime graph $\Gamma(G)$?
    \item  For which groups are there only finitely many groups with the same prime
    graph $\Gamma(G)$ as $G$?
    \item  Which groups $G$ are uniquely determined by isomorphism type of their prime graph?
\end{enumerate}

Similar to the previous discussion about the function $h$ (See Definition \ref{2.3}). In \cite{29} the authors gave the following definition (see \cite{29}, page 188).

\begin{defi}
We say that the group $G$ is
\begin{itemize}
\item recognizable \textup{(}by prime graph\textup{)} if for every group $H$ the equality
$\Gamma(H) =\Gamma(G)$ implies that $G\cong H$;
\item $k$-recognizable \textup{(}by prime graph\textup{)}, where $k$ is a positive integer, if there are exactly $k$ pairwise non-isomorphic groups having the same prime graph as $G$;
\item  almost recognizable \textup{(}by prime graph\textup{)} if it is $k$-recognizable by prime graph for a positive integer $k$;
\item  unrecognizable \textup{(}by prime graph\textup{)} if there are infinitely many pairwise non-isomorphic groups having the same prime graph as $G$.
\end{itemize}
\end{defi}

Next, the following theorem can be found in \cite{MPS}.

\begin{thm}
Every finite simple exceptional group of Lie type, which is isomorphic
to neither ${}^2B_2(2^{2n+1})$ with $n\geqslant 1$ nor $G_2(3)$ and whose prime graph has at least three connected components, is almost recognizable by prime
graph. Moreover, groups ${}^2B_2(2^{2n+1})$, where $n\geqslant 1$, and $G_2(3)$ are unrecognizable by prime graph.
\end{thm}

There are differences between having the same prime graph and having an isomorphic prime graph. For example, $\Gamma(A_{10})\neq\Gamma(\Aut(J_2))$ but $\Gamma(A_{10})$ and $\Gamma(\Aut(J_2))$ are isomorphic as abstract graphs.

\begin{pro}\label{6.17}
 Which simple groups can be characterized by their prime graphs and orders?
\end{pro}

The above question is a generalization of the work that characterizes all simple groups using ``two orders". Consider the classical groups $B_n(q)$ and $C_n(q)$ (See \cite{176}), when $n\geqslant 3$ and $q$ is odd, their prime graphs and orders are same. Therefore, not all simple groups can be characterized by their prime graphs and orders.

\subsection{Graphs defined on groups}\

There are  various graphs defined on a group, but the one most closely associated with quantitative characterization is the prime graph (GK graph). Cameron \cite{27} introduced the following graphs: Their vertex sets are finite groups $G$, and their edges reflect the structure of the groups in a certain way, such that the automorphism group of $G$ related to the automorphism of the graph. These graphs include the commuting graph (two vertices $x$ and $y$ are joined if and only if $xy=yx$;
this graph was first introduced in 1955, see \cite{22}),
the generating graph ($x$ and $y$ joined if and only if $\langle x, y\rangle=G$; the generating graph was first introduced in 1996, see reference \cite{23,109}),
the power graph (vertices $x$ and $y$ are joined if and only if one of $x$ and $y$ is a power of the other; this graph was first studied in 2000, see \cite{93,139,KSC}),
the enhanced power graph (two vertices $x$ and $y$ are joined if and only if $\langle x, y\rangle$ is not cyclic; this graph was first studied in 2007, see \cite{2, 3})
and the deep commuting graph (two vertices, $x$ and $y$, are joined if and only if their preimages, $\overline{x}$ and $\overline{y}$, commute in every central extension of $G$. In other words, for every group $H$ with a central subgroup $Z$ such that $H/Z\cong G$, and their preimages $\overline{x}$ and $\overline{y}$ commute in $H$, see \cite{28}).

This paper mainly discusses the quantitative characterization of finite simple groups, which is an essential component of quantitative properties of groups and a significant topic of group theory research. Studying finite groups using the set of element orders and characterizing finite simple groups by ``two orders" were first proposed and investigated by the author, which have a significant influence on the field of group theory at home and abroad.
Since the 1980s, the author has successively worked at Southwest Normal University (including a concurrent doctoral supervisor position at Sichuan University), Suzhou University, and Chongqing University of Arts and Sciences. Additionally, the author has received funding for 11 projects from the National Natural Science Foundation of China, enabling continuous progress in this work.
Our intention in writing this article is not only to systematically summarize our work, distill new viewpoints, ideas, and methods to solve related problems, but also to further delve into this work and find more significant applications.

\section*{Acknowledgement}

We sincerely thank the reviewers for their helpful suggestions for revisions. Ping Jin and Guohua Qian provided valuable advice on writing and specific content for this paper. Nanying Yang and Rulin Shen also offered their revision suggestions for this paper. In particular, Jinbao Li dedicated a considerable amount of effort to the editing and revising of this paper. The author would like to extend his heartfelt gratitude to all of them for their contributions. Finally, the author would like to express his gratitude to Andrey Vasil'ev, Mariya Grechkoseeva and Natalia Maslova for their valuable suggestions on the revised version. The author thanks Xiaofang Gao for translating this article from Chinese version to revised version, which is helpful for the author's communication at home and abroad.

The author is very grateful to the referees for providing valuable revision suggestions and specific recommendations.

\end{document}